\documentclass[final,oneside, 12pt, a4paper]{amsart}
\usepackage{todonotes}
\PassOptionsToPackage{final}{graphix} 

\RequirePackage{xparse}     

\RequirePackage[utf8]{inputenc} 
\RequirePackage[T1]{fontenc}    
\RequirePackage[a4paper,margin=3cm]{geometry}
\usepackage[sb]{libertinus}
\usepackage[amsthm,lcgreekalpha]{libertinust1math}
\usepackage[cal=boondoxo, calscaled=0.96, frak=esstix, frakscaled=0.93]{mathalfa}

\RequirePackage[tracking=true]{microtype} 
\SetTracking   
{encoding=T1, shape=sc}
{20}
\SetProtrusion 
{encoding=T1,size={7,8}}
{1={ ,750},2={ ,500},3={ ,500},4={ ,500},5={ ,500},
  6={ ,500},7={ ,600},8={ ,500},9={ ,500},0={ ,500}}

\RequirePackage{amssymb}  
\RequirePackage{mathrsfs} 
\RequirePackage{mathbbol} 

\RequirePackage{mathtools}  
\RequirePackage{leftidx}  
\RequirePackage{bussproofs} 
\RequirePackage{scalerel}   
\RequirePackage{faktor}     
\RequirePackage{xfrac}      

\usepackage{ytableau}
\ytableausetup{aligntableaux=top,mathmode}

\DeclareFontFamily{U}{DSSerif}{\skewchar \font =45}
\DeclareFontShape{U}{DSSerif}{m}{n}{<-> s*[1]  DSSerif}{}
\DeclareFontSubstitution{U}{DSSerif}{m}{n}
\DeclareMathAlphabet{\mathbbbb}{U}{DSSerif}{m}{n}

\RequirePackage[english]{babel} 
\RequirePackage{csquotes}       

\RequirePackage{hyphenat} 

\definecolor{red}{RGB}{175, 49, 39}             

\RequirePackage{xspace}         
\RequirePackage{enumitem}
\RequirePackage{aligned-overset} 
\usepackage[style=alphabetic,natbib=true, maxbibnames=2]{biblatex}
\DeclareSourcemap{
  \maps[datatype=bibtex]{
    \map{
      \step[fieldset=doi, null]
      \step[fieldset=url, null]
      \step[fieldset=urldate, null]
      \step[fieldset=language, null]
      \step[fieldset=primaryClass, null]
      \step[fieldset=note, null]
      \step[fieldset=pages, null]
      \step[fieldset=issn, null]
      \step[fieldset=isbn, null]
      \step[fieldset=school, null]
    }
  }
}
\DefineBibliographyStrings{english}{%
  andothers = {et\,al\adddot},
  and       = {\addcomma\space},
  in = {},
}

\RequirePackage[colorlinks,          
linktoc=page,
linkcolor={red},
citecolor={red},
urlcolor={red}
]{hyperref}
\RequirePackage{amsthm}

\newlist{thmlist}{enumerate}{1}
\setlist[thmlist]{label=(\roman{thmlisti}), ref=\thetheorem.(\roman{thmlisti}), noitemsep}
\newlist{thmlist*}{enumerate}{1}
\setlist[thmlist*]{label=(\roman{thmlist*i}), ref=(\roman{thmlist*i}), noitemsep}

\numberwithin{equation}{section}

\theoremstyle{plain}
\newtheorem{theorem}{Theorem}
\newtheorem*{theorem*}{Theorem}
\newtheorem{lemma}[theorem]{Lemma}
\newtheorem{corollary}[theorem]{Corollary}
\newtheorem{proposition}[theorem]{Proposition}
\newtheorem*{proposition*}{Proposition}

\theoremstyle{definition}
\newtheorem{definition}[theorem]{Definition}
\newtheorem{example}[theorem]{Example}

\newtheorem{remark}[theorem]{Remark}
\newtheorem{construction}[theorem]{Construction}
\newtheorem{notation}[theorem]{Notation}

\newtheorem*{question*}{Question}
\newtheorem*{remark*}{Remark}

\newcommand*{\eg}{e.g.\xspace}

\renewcommand*{\colon}{\!\nobreak\mskip2mu\mathpunct{}\nonscript%
  \mkern-\thinmuskip{:}\mskip6muplus1mu\relax%
}
\newcommand*{\from}{\ensuremath{\colon}} 
\renewcommand*{\to}{\ensuremath{\longrightarrow}}
\renewcommand*{\mapsto}{\ensuremath{\longmapsto}}

\DeclareMathOperator{\op}{op}

\makeatletter
\newcommand{\lMod}[1]{\ensuremath{{\mathsf{Mod}}_{#1}}}

\DeclareMathOperator{\Hom}{Hom}      
\DeclareMathOperator{\Aut}{Aut}      

\DeclareMathOperator{\Res}{Res}
\DeclareMathOperator{\Ind}{Ind}

\newcommand{\Rep}[1]{\ensuremath{\mathsf{Rep}_{#1}}}

\let\oldsum\sum
\renewcommand{\sum}{\oldsum\nolimits}

\author{Sebastian Halbig}
\address{Sebastian Halbig, Philipps-Universität Marburg, Arbeitsgruppe Algebraische Lie-Theorie, Hans-Meerwein-Strasse 6, 35043 Marburg, Germany}
\email{Sebastian.Halbig@uni-marburg.de}

\author{Christian Lomp}
\address{Christian Lomp, CMUP, Departamento de Matemática, Faculdade de Ciências,
  Universidade do Porto, Rua do Campo Alegre s/n, 4169–007 Porto,
  Portugal.}
\email{clomp@fc.up.pt}
\thanks{This work was initiated during a visit of S.H. at the Universidade do Porto in October 2025. S.H. thanks Paula Carvalho and Christian Lomp for their warm hospitality.
  The second author was partially supported by CMUP – Centro de Matemática da Universidade do Porto,
  member of LASI, which is financed by national funds through FCT – Fundação para a
  Ciência e a Tecnologia, I.P., under the project with reference UID/00144/2025, doi: https://doi.org/10.54499/UID/00144/2025}

\title{Irreducible representations of generalised Kac--Paljutkin Hopf algebras}
\date{\today}

\pdfoutput=1
\begin{document}
\begin{abstract}
  The aim of this note is to provide a self-contained classification of the irreducible representations of  generalised Kac--Paljutkin Hopf algebras, recently introduced by the second author.
\end{abstract}
\maketitle

\section{Introduction}\label{sec:introduction}
The Kac--Paljutkin Hopf algebra, introduced in \cite{kats-palyutkin1966:FiniteRingGroups},
is the dimension-wise smallest complex non-commutative and non-cocommutative Hopf algebra which is semisimple and cosemisimple.
This construction was extended by Pansera to a class of \(2n^{2}\)-dimensional Hopf algebras \cite{pansera2019:HopfAlgebrasQuantumPolynomial} whose irreducible representations were considered in~\cite{chen-yang-wang2021:GrothendieckHopfKacPaljutkin}.
Recently, the second author provided in~\cite{lomp2025:GeneralizedKacPaljutkin} a further generalisation, leading to a family  \((H_{n,m})_{n,m\geq 2}\)  of \(n^{m}m!\)-dimensional non-commutative, non-cocommutative, semisimple, and cosemisimple Hopf algebras.
Unlike Pansera's examples where all irreducible representations are either one or two-dimensional, the situation for the  generalised Kac--Paljutkin Hopf algebras \(H_{n,m}\) is more involved.
For example, the representations \(\Rep{S_{m}}\) of the symmetric group \(S_{m}\) embed into \(\lMod{H_{n,m}}\), see Remark~\ref{rmk:embedding-of-sym-grp-reps}.
In this short article, we are therefore going to classify the isomorphism classes of the  irreducible representations of \(H_{n,m}\) by group theoretic means.

\begin{theorem}[Theorem~\ref{thm:the-one-we-want}]\label{thm:main-result}
  Let \(  \Bbbk\) be an algebraically closed field of characteristic zero and fix two positive integers \(2 \leq n,m \in \mathbb{N}\).
  The generalised Kac--Paljutkin \( \Bbbk\)-Hopf algebra \(H_{n,m}\) is isomorphic as an algebra to the group algebra \( \Bbbk[ \mathbb{Z}_{n} \wr S_{m}]\).

  Its irreducible representations are in bijection with \(n\)-tuples of partitions
  \begin{equation}
    (\lambda_{1} \vdash k_{1}, \lambda_{2} \vdash k_{2}, \dots , \lambda_{n} \vdash k_{n}) \quad \text{such that } k_{1}+ \dots + k_{n}=m.
  \end{equation}
\end{theorem}

The article is structured as follows.
In Section~\ref{sec:KP-algebras}, we recall the definition of the generalised Kac--Paljutkin Hopf algebras and prove that they are isomorphic, as algebras (but not Hopf algebras), to certain group algebras, see Proposition~\ref{prop:KP-is-hyperoct}.
To give a self-contained proof of our main theorem, we briefly recall the necessary parts of the representation theory of finite groups and how to determine irreducible representations of semidirect products using Clifford theory in Section~\ref{sec:wreath-products}.
Finally, we prove our main result,  Theorem~\ref{thm:the-one-we-want}, in Section~\ref{sec:proof-of-the-main-theorem}.
\numberwithin{theorem}{section}

\section{Generalised Kac--Paljutkin algebras and wreath products}\label{sec:KP-algebras}
Throughout this note, we work over an algebraically closed field \( \Bbbk\) of characteristic zero and fix natural numbers \(2 \leq n, m \in \mathbb{N}\) as well as a primitive \(2n\)-th root of unity \(\zeta \in  \Bbbk\).
For convenience, we write \(q \eqdef \zeta^{2}\) and \([k] \eqdef \{0, … , k-1\}\) for all \(k\in \mathbb{N}\).
Moreover, given \(1\leq l \leq m-1\), we set  \(\sigma_{l} \eqdef (l \; l+1)\in S_{m}\).
\medskip

We define the generalised Kac--Paljutkin Hopf algebras following  \cite[Theorem~7]{lomp2025:GeneralizedKacPaljutkin}.
\begin{definition}\label{def:GKP}
  The \emph{generalised Kac--Paljutkin Hopf algebra} \(H_{n,m}\) is generated by elements
  \(x_{1}, \dots , x_{m}, z_{1}, \dots , z_{m-1}\) satisfying for every \(1\leq i,j \leq m\) and \(1\leq k, l \leq m-1\) the relations
  \begin{subequations}
    \begin{gather}
      x_{i}^{n}=1, \qquad
      x_{i} x_{j} = x_{j} x_{i}, \qquad
      z_{l} x_{i} = x_{\sigma_{l}(i)} z_{l}, \label{eq:KP-rel-1}\\
      z_{l}z_{k}= z_{k} z_{l}\;\;   \text{if }|k-l| \geq 2, \quad z_{l}z_{l+1}z_{l}= z_{l+1}z_{l}z_{l+1} \;\; \text{ for } 1\leq l\leq m-2,\\
      z_{l}^{2} = \frac{1}{n} \sum_{i,j=1}^{n-1} q^{-ij} x^{i}_{l}x^{j}_{l+1}.
    \end{gather}
  \end{subequations}
  Furthermore, the comultiplication and antipode of \(H_{n,m}\) are given by
  \begin{subequations}
    \begin{gather}
      \Delta(x_{i}) = x_{i} \otimes x_{i}, \qquad \Delta(z_{l}) =\bigg(\frac{1}{n}\sum_{i,j=0}^{n-1} q^{-ij} x_{l}^{i} \otimes x_{l+1}^{j}\bigg) (z_{l} \otimes z_{l}), \\
      S(x_{i}) = x_{i}^{n-1}, \qquad
      S(z_{l}) = z_{l}.
    \end{gather}
  \end{subequations}
\end{definition}
While the generators \(z_{1}, \dots , z_{m-1}\) do not generate a subgroup isomorphic to \(S_{m}\) in \(H_{n,m}\), they can be used to embed the representations \(\Rep{S_{m}}\) of \(S_{m}\) into the category \(\lMod{H_{n,m}}\) of \(H_{n,m}\)-modules.
\begin{remark}\label{rmk:embedding-of-sym-grp-reps}
  We can form in \(H_{n,m}\) the two-sided ideal \(I\) spanned by \(\{x_{i}-1 \mid 1\leq i \leq m\}\).
  A direct calculation shows that it is a Hopf ideal and the quotient \(H_{n,m}/I\) is isomorphic to the group algebra \(  \Bbbk[S_{m}]\).
  Therefore, the canonical projection \(\pi\from H_{n,m}\to H_{n,m}/I \cong  \Bbbk[S_{m}]\) induces a monoidal embedding
  \begin{equation*}
    \pi^{*} \from \Rep{S_{m}} \to \lMod{H_{n,m}}, \qquad\quad (V,\rho) \mapsto (V, \rho \pi).
  \end{equation*}

\end{remark}

The elements \(x_{1}, \dots x_{m}\) of \(H_{n,m}\) span a subgroup isomorphic to the \(m\)-fold Cartesian product of the cyclic group \( \mathbb{Z}_{n} \eqdef \mathbb{Z} / n \mathbb{Z}\), which we will simply denote by \( \mathbb{Z}_{n}^{m}\).
In order to classify the irreducible representations of \(H_{n,m}\), we first need to study the representation theory of \(\mathbb{Z}_{n}^{m}\).
Given any \(m\)-tuple \(\mathbf{\lambda} \eqdef(\lambda_{1}, \dots , \lambda_{m}) \in \mathbb{Z}_{n}^{m}\), there is a one-dimensional representation \( \Bbbk_{\mathbf{\lambda}}\) defined by \(x_{i} \cdot 1_{  \Bbbk} = q^{-\lambda_{i}}\) and any irreducible representation of \(\mathbb{Z}_{n}^{m}\) is isomorphic to some \(  \Bbbk_{\mathbf{\lambda}}\).
The next Lemma is a standard result in the representation theory of finite groups, see \eg~\cite{fulton-harris1991:RepresentationTheory, isaacs2006:CharacterTheoryFiniteGroups,schedler2021:GroupRepresentationTheory}.
It concerns the idempotents of \( \Bbbk[\mathbb{Z}_{n}^{m}]\) corresponding to these representations.
For convenience,  we will use for all \(\mathbf{\lambda}, \mathbf{i} \in \mathbb{Z}_{n}^{m}\) the shorthand notations
\begin{equation*}
  \mathbf{\lambda} \cdot \mathbf{i} \eqdef \lambda_{1}i_{1}+ \dots+ \lambda_{m}i_{m}\qquad \text{and} \qquad
  \mathbf{x}^{\mathbf{i}} \eqdef x_{1}^{i_{1}} x_{2}^{i_{2}} \cdot \dots \cdot x_{m}^{i_{m}} \in H_{n,m}.
\end{equation*}

\begin{lemma}\label{lem:irreducible-reps-of-A}
  Given \(\mathbf{\lambda}\in \mathbb{Z}_{n}^{m}\), define
  \begin{equation}
    \Lambda_{\mathbf{\lambda}} \eqdef \frac{1}{n^{m}} \sum_{\mathbf{i} \in \mathbb{Z}_{n}^{m}}q^{\mathbf{\lambda}\cdot \mathbf{i}}\mathbf{x}^{\mathbf{i}} \in H_{n,m}
  \end{equation}
  For each \(1\leq i\leq m\), we have  \(x_{i} \Lambda_{\mathbf{\lambda}} = q^{-\lambda_{i}} \Lambda_{ \mathbf{\lambda}}\) and the set \(\{\Lambda_{\mathbf{\lambda}} \mid \mathbf{\lambda} \in \mathbb{Z}_{n}^{m}\}\) is a basis of \(  \Bbbk[\mathbb{Z}_{n}^{m}] \subset H_{n,m}\) consisting of primitive (central) orthogonal idempotents.
\end{lemma}
\begin{definition}\label{def:AW-basis}
  We call \(\{\Lambda_{\mathbf{\lambda}} \mid \mathbf{\lambda} \in \mathbb{Z}_{n}^{m}\}\) the \emph{Artin--Wedderburn basis} of \(  \Bbbk[\mathbb{Z}_{n}^{m}]\).
\end{definition}
The next result is a generalisation of Remark 2.14 of \cite{masuoka1995:SemisimpleHopf}.
\begin{lemma}\label{lem:sq-root-wong-relation}
  For each \(1\leq l \leq m-1\) set
  \begin{equation}\label{eq:the-element-that-fixes-it}
    y_{l} = \sum_{\mathbf{\lambda} \in \mathbb{Z}_{n}^{m}} \zeta^{-\lambda_{l}\lambda_{l+1}} \Lambda_{\lambda} \qquad \text{ and }\qquad
    s_{l} = y_{l}z_{l}.
  \end{equation}
  \begin{thmlist}
    \item The elements  \((y_{l})_{1\leq l \leq m-1}\) are units of \(H_{n,m}\) of order \(2n\), and
    \item for all \(1\leq i \leq m\) and \(1\leq k,l\leq m-1\) we have
    \begin{subequations}
      \begin{gather}
        s_{l} x_{i} = x_{\sigma_{l}(i)}s_{l}, \qquad
        s_{l}^{2} = 1, \label{eq:sym-grp-1} \\
        s_{l} s_{k} = s_{k} s_{l} \;\; \text{if }|k-l|\geq 2, \quad
        s_{l} s_{l+1} s_{l} = s_{l+1} s_{l} s_{l+1} \;\;\text{if } l \leq m-2. \label{eq:sym-grp-2}
      \end{gather}
    \end{subequations}
  \end{thmlist}
\end{lemma}
\begin{proof}
  Using  the orthogonality of the $\Lambda_\lambda$'s,  we have  $y_{l}^k = \sum_{\mathbf{\lambda} \in \mathbb{Z}_{n}^{m}} \zeta^{-k\lambda_{l}\lambda_{l+1}} \Lambda_{\lambda}$, for all $k\geq0$.
  In particular, as $\zeta$ has order $2n$, we have for all \(1\leq l \leq m-1\) and \(k\geq 0\) that
  \begin{equation*}
    y_{l}^{k}= \sum_{\mathbf{\lambda} \in \mathbb{Z}_{n}^{m}} \zeta^{-k \lambda_{l}\lambda_{l+1}} \Lambda_{\mathbf{\lambda}} =  \sum_{\mathbf{\lambda} \in \mathbb{Z}_{n}^{m}} \Lambda_{\mathbf{\lambda}} = 1 \iff  2n | k.
  \end{equation*}

  Before we can prove Equations~\eqref{eq:sym-grp-1} and~\eqref{eq:sym-grp-2}, we need to establish two auxiliary identities.

  First, as \(z_{l}^{2} \in  \Bbbk[ \mathbb{Z}_{n}^{m}]\), we can express it in terms of the Artin--Wedderburn basis of \(  \Bbbk[ \mathbb{Z}_{n}^{m}]\).
  To that end, we note that for any \(\mathbf{\lambda} \in \mathbb{Z}_{n}^{m}\), we have
  \begin{align*}
    z_{l}^{2} \Lambda_{\mathbf{\lambda}}
    & = \frac{1}{n}\sum_{i,j=0}^{n-1} q^{-ij}x_{l}^{i}x_{l+1}^{j} \Lambda_{\mathbf{\lambda}}
      = \frac{1}{n} \sum_{i,j=0}^{n-1}q^{-ij - \lambda_{l}i - \lambda_{l+1}j} \Lambda_{\mathbf{\lambda}}
      = q^{ \lambda_{l} \lambda_{l+1}} \Lambda_{\mathbf{\lambda}},
  \end{align*}
  because we have for fixed \(a, b \in \mathbb{Z}_{n}\) :
  \begin{align*}
    \sum_{i,j=0}^{n-1} q^{-ij -ai - bj}
    & = \sum_{i,j=0}^{n-1} q^{-(i+b)(j+a) +ab}
      = q^{ab} \sum_{r,s=0}^{n-1}q^{-rs} 
      = n q^{ab}.
  \end{align*}
  Thus we get
  \begin{equation} \label{eq:the-equation-that-explains-it-all}
    z_{l}^{2}
    = \sum_{\mathbf{\lambda}\in \mathbb{Z}_{n}^{m}}  z_{l}^2\Lambda_{\mathbf{\lambda}}
    = \sum_{\mathbf{\lambda}\in \mathbb{Z}_{n}^{m}} q^{\lambda_{l} \lambda_{l+1}} \Lambda_{\mathbf{\lambda}}
    = \sum_{\mathbf{\lambda}\in \mathbb{Z}_{n}^{m}} \zeta^{2\lambda_{l} \lambda_{l+1}} \Lambda_{\mathbf{\lambda}}
    = y^{-2}_{l}.
  \end{equation}

  Second, given any \(1\leq l \leq m-1\) as well as  \(\mathbf{i}\in \mathbb{Z}_{n}^{m}\) and using the short-hand notation \(\sigma_{l}(\mathbf{i})\eqdef (i_{\sigma_{l}(1)}, i_{\sigma_{l}(2)}, \dots,i_{\sigma_{l}(m)})\), we compute
  \begin{equation}\label{eq:z-action-on-simples}
    z_{l} \Lambda_{\mathbf{\lambda}}
    = \frac{1}{n^{m}}\sum_{\mathbf{i} \in \mathbb{Z}_{n}^{m}} q^{\mathbf{\lambda}\cdot \mathbf{i}}\mathbf{x}^{\sigma_{l}(\mathbf{i})}z_{l}
    = \frac{1}{n^{m}} \sum_{\mathbf{i} \in \mathbb{Z}_{n}^{m}} q^{\sigma_{l}(\mathbf{\lambda})\cdot \mathbf{i}}\mathbf{x}^{\mathbf{i}} z_{l}
    = \Lambda_{\sigma_{l}(\mathbf{\lambda})} z_{l},
  \end{equation}
  implying for all \(1\leq k,l \leq m-1\) that
  \begin{equation}\label{eq:how-z-acts-on-y}
    z_{l} y_{k}
    = \sum_{\mathbf{\lambda}\in \mathbb{Z}_{n}^{m}} \zeta^{-\lambda_{k}\lambda_{k+1}} \Lambda_{\sigma_{l}(\mathbf{\lambda})} z_{l}
    = \sum_{\mathbf{\lambda}\in \mathbb{Z}_{n}^{m}} \zeta^{-\lambda_{\sigma_{l}(k)}\lambda_{\sigma_{l}(k+1)}} \Lambda_{\mathbf{\lambda}} z_{l}.
  \end{equation}

  With these identities established, we now prove that Equations~\eqref{eq:sym-grp-1} and~\eqref{eq:sym-grp-2} hold:
  \begin{gather*}
    s_{l} x_{i}
    = y_{l} z_{l} x_{i}
    = y_{l} x_{\sigma(l)} z_{i}
    =  x_{\sigma(l)} y_{l} z_{i}
    = x_{\sigma(l)} s_{l}, \\
    s_{l}^{2}
    = y_{l} z_{l} y_{l} z_{l}
    \overset{{\eqref{eq:how-z-acts-on-y}}}
    = y_{l}^{2} z_{l}^{2}
    \overset{{\eqref{eq:the-equation-that-explains-it-all}}}
    = 1, \\
    s_{l}s_{k}
    = y_{l} z_{l} y_{k} z_{k}
    \overset{{\eqref{eq:how-z-acts-on-y}}}
    = y_{k} z_{k} y_{l} z_{l}
    = s_{k} s_{l}, \text{ if } |k-l|\geq 2.
  \end{gather*}
  It remains to show \(s_{k}s_{k+1}s_{k}= s_{k+1} s_{k} s_{k+1}\).
  Iteratively appling Equation~\eqref{eq:how-z-acts-on-y}, leads to
  \begin{align*}
    & s_{k} s_{k+1}s_{k} = y_{k} z_{k} y_{k+1}z_{k+1} y_{k} z_{k} \\
    & = \big( \sum_{\mathbf{\lambda} \in \mathbb{Z}_{n}^{m}} \zeta^{-\lambda_{k}\lambda_{k+1}} \Lambda_{\mathbf{\lambda}} \big) \big( \sum_{\mathbf{\mu} \in \mathbb{Z}_{n}^{m}} \zeta^{-\mu_{k}\mu_{k+2}} \Lambda_{\mathbf{\mu}} \big) \big( \sum_{\mathbf{\tau} \in \mathbb{Z}_{n}^{m}} \zeta^{-\tau_{k+1}\tau_{k+2}} \Lambda_{\mathbf{\tau}} \big)z_{k}z_{k+1}z_{k}\\
    & s_{k+1}s_{k}s_{k+1} = y_{k+1} z_{k+1} y_{k}z_{k} y_{k+1} z_{k+1} \\
    & = \big( \sum_{\mathbf{\lambda} \in \mathbb{Z}_{n}^{m}} \zeta^{-\lambda_{k+1}\lambda_{k+2}} \Lambda_{\mathbf{\lambda}} \big) \big( \sum_{\mathbf{\mu} \in \mathbb{Z}_{n}^{m}} \zeta^{-\mu_{k}\mu_{k+2}} \Lambda_{\mathbf{\mu}} \big) \big( \sum_{\mathbf{\tau} \in \mathbb{Z}_{n}^{m}} \zeta^{-\tau_{k}\tau_{k+1}} \Lambda_{\mathbf{\tau}} \big) z_{k+1}z_{k}z_{k+1} \\
  \end{align*}
  Thus, \(s_{k}s_{k+1}s_{k}=s_{k+1}s_{k}s_{k+1}\) follows from \(z_{k} z_{k+1} z_{k} = z_{k+1} z_{k} z_{k+1}\).
\end{proof}
We will now use the set \(\{x_{1}, \dots , x_{m}, s_{1}, \dots , s_{m-1}\}\) to derive a presentation of \(H_{n,m}\) in terms of a group algebra.

\begin{definition}\label{def:wreath-product}
  Let \(A\) be a group and \(G \subset S_{m}\) a subgroup of a symmetric group.
  The \emph{wreath product} \(A \wr G\) is the group with underlying set \(A^{m} \times G\) and the multiplication
  \begin{equation}
    (a_{1}, \dots , a_{m},g)(b_{1}, \dots , b_{m}, h) =
    (a_{1}b_{g^{-1}(1)}, \dots , a_{m} b_{g^{-1}(m)}, gh).
  \end{equation}
\end{definition}
The terminology of the next definition follows \cite{osima1954:RepresentationGeneralisedSymmetricGroup}.
\begin{definition}\label{def:generalised-hyoct-group}
  We refer to \(\mathbb{Z}_{n} \wr S_{m}\) as the \emph{generalised symmetric group}.
\end{definition}

In the Sheppard--Todd classification of complex reflection groups, the group  \( \mathbb{Z}_{n} \wr S_{m}\) is denoted by \(G(n,1,m)\), see \cite{shephard-todd1954:FiniteUnitaryReflectionGroups}.
The groups  \( \mathbb{Z}_{2} \wr S_{m}\) are also called \emph{signed symmetric groups} or \emph{hyperoctahedral groups}.

Since wreath products are special semidirect products, we can describe them in terms of generators and relations.

\begin{remark}\label{rmk:presentation-of-hyperoct-grp}
  The generalised symmetric group \(\mathbb{Z}_{n} \wr S_{m}\) has a presentation in terms of the generators \(a_{1}, \dots , a_{m}, b_{1}, \dots , b_{m-1}\) and relations given for all \(1\leq i,j \leq m\) and \(1\leq k, l \leq m-1\) by
  \begin{gather*}
    a_{i}^{n}=1 \qquad
    a_{i} a_{j} = a_{j} a_{i} \qquad
    b_{l} a_{i} = a_{\sigma_{l}(i)} b_{l}, \\
    b_{l}^{2} =1, \qquad
    b_{l}b_{k} = b_{k} b_{l} \;\; \text{if } |k-l|\geq2, \quad
    b_{l}b_{l+1}b_{l} = b_{l+1} b_{l} b_{l+1} \;\; \text{if } l\leq m-2.
  \end{gather*}
\end{remark}
Combining Remark~\ref{rmk:presentation-of-hyperoct-grp} with Lemma~\ref{lem:sq-root-wong-relation}, leads to an identification of \(H_{m,n}\) with \( \Bbbk[\mathbb{Z}_{n}\wr S_{m}]\).
\begin{proposition}\label{prop:KP-is-hyperoct}
  There exists an isomorphism of algebras \(\phi \from \Bbbk[\mathbb{Z}_{n}\wr S_{m}] \to H_{n,m}\) satisfying
  \begin{equation*}
    \phi(a_{i})= x_{i}, \qquad \phi(b_{l})=s_{l} \qquad \text{for all } 1\leq i \leq m \text{ and } 1\leq l \leq m-1.
  \end{equation*}
\end{proposition}

Note that the previous result only concerns the algebra but not the coalgebra structure.

\begin{remark}\label{rmk:not-Hopf algebra isomorphism}
  If \(\Bbbk[\mathbb{Z}_{n}\wr S_{m}]\) and \(H_{n,m}\) were isomorphic as Hopf algebras, it would imply the contradictory claim that \(\Delta(z_{l})= \Delta^{\op}(z_{l})\) for all \(1\leq l \leq m-1\).
\end{remark}

Isomorphism classes of irreducible representations of wreath products have been extensively studied using various setups, see for example \cite{james-kerber1981:RepresentationTheorySymmetricGroup, mishra-srinivasan2016:OkounkovVershikRepresentationTheoryWreathProducts, paul-pfeiffer2025:ComputingYoungRepresentationsGeneralisedSymmetricGroups,ceccherini-silberstein-scarabotti-tolli2022:RepresentationFiniteGroupExtensions}.
In order to keep our exposition self-contained, we will provide in the next section a brief overview of the representation theory of wreath products and state an elementary proof of our main theorem in Section~\ref{sec:proof-of-the-main-theorem}.

\section{Representations of wreath products}\label{sec:wreath-products}

To determine the irreducible representations  of \(H_{n,m}\), we briefly recall certain standard group representation-theoretic tools.
Our exposition follows the books~\cite{ceccherini-silberstein-scarabotti-tolli2022:RepresentationFiniteGroupExtensions,howe2022:InvitationToRepresentationTheory}.

\subsection{Clifford theory}\label{sec:clifford-theory}
Throughout this section, let \(G\) be a finite group, \(H\) a subgroup, and \(N\) a normal subgroup of \(G\), unless specified explicitly otherwise.

\begin{notation}\label{not:simple-representations}
  Given a group \(G\), we write \(\Rep{G}\) for its category of representations and \(\widehat{G}\) for the set of isomorphism classes of irreducible representations.
\end{notation}

\begin{definition}\label{def:reduced-induced-reps}
  Suppose \(H\) is a subgroup of \(G\) and  consider representations \(V\in \Rep{H}\) as well as \(W\in \Rep{G}\).
  We call
  \begin{equation}
    \Ind_{H}^{G}V \eqdef  \Bbbk [G] \otimes_{ \Bbbk[H]} V \in \Rep{G} \quad \text{and} \quad
    \Res^{G}_{H}W \eqdef \Hom_{ \Bbbk[G]}( \Bbbk[G], W) \in \Rep{H}
  \end{equation}
  the \emph{induction of \(V\) from \(H\) to \(G\)} and the \emph{restriction of \(W\) from \(G\) to \(H\)}, respectively.
\end{definition}

Clifford theory studies the behaviour of induction and restriction between \(G\) and normal subgroups \(N\) of \(G\).

\begin{definition}\label{def:conjugation-action-on-reps}
  For any \(g\in G\) and any \((V,\rho) \in \Rep{N}\), we write \({}^{g}V \in \Rep{N}\) for the representation defined by
  \begin{equation*}
    ^{{g}}\rho \from N \to \Aut_{ \Bbbk}(V), \qquad n \mapsto\rho(g^{-1}ng).
  \end{equation*}
\end{definition}
Note that this gives rise to an action of \(G\) on \(\widehat{N}\).

\begin{definition}\label{def:inertia-subgroup}
  We call the stabiliser \(I_{G}([V])\) of \([V] \in \widehat{N}\) under the action of \(G\) the \emph{inertia subgroup} of \([V]\).
\end{definition}

The restriction \(\Res^{G}_{N} W\) of an irreducible representation \(W\) of \(G\) has  a unique decomposition into \emph{isotypic components}
\begin{equation*}
  \Res^{G}_{N} W \cong \bigoplus_{[V]\in \widehat{N}} U_{[V]}, \qquad \text{where }
  U_{[V]}\cong V^{a_{[V]}}.
\end{equation*}
\begin{definition}\label{def:multiplicity}

  Given a representation \(W\) of \(G\) and an irreducible representation \(V\) of \(N\), we write
  \begin{equation*}
    [V:W]_{N} \eqdef \dim \Hom_{N}(V, \Res^{G}_{N} W)
  \end{equation*}
  for the \emph{multiplicity} of \(V\) in \(\Res^{G}_{N}W\).
    We denote by \(\widehat{G}_{[V]}\) the set of all isomorphism classes of irreducible representations \(W\) of \(G\) such that \([V:W]_{N}\geq 1\).
\end{definition}

The next result is a variant of Clifford's theorem.
A proof can be found for example in \cite[Theorems~2.5(3) and  2.12]{ceccherini-silberstein-scarabotti-tolli2022:RepresentationFiniteGroupExtensions}.
\begin{lemma}\label{lem:Cliffords-thm}
  Let \(V\in \Rep{N}\) and \(W\in \Rep{G}\) be irreducible representations of \(N\) and \(G\), respectively.
  If \(a= [V:W]_{N}\geq 1\), we have
  \begin{equation}
    \Res^{G}_{N} W  \cong \bigoplus_{[h] \in G/I_{G}([V])} {}^{h}(V^{a}).
  \end{equation}
  Moreover, there is a \(U \in \Rep{I_{G}([V])}\) such that \(\Res^{I_{G}([V])}_{N} U\cong V^{a}\) and \(W \cong \Ind_{I_{G}([V])}^{G}U\).
\end{lemma}

Subsequently, we can obtain \(\widehat{G}\) by studying the representation theory of inertia subgroups, leading to the \emph{Clifford correspondence}, see \cite[Theorem~2.12]{ceccherini-silberstein-scarabotti-tolli2022:RepresentationFiniteGroupExtensions}.
\begin{lemma}\label{lem:clifford-correspondence}
  Fix an irreducible representation \(V\) of \(N\).
  The map
  \begin{equation}
    \widehat{I_{G}([V])}_{[V]} \to \widehat{G}_{[V]}, \qquad\quad
    [U] \mapsto [\Ind_{I_{G}([V])}^{G} U]
  \end{equation}
  is a bijection.
\end{lemma}

As wreath products are special cases of semidirect products, their irreducible representations can be described using the \emph{little groups method}, see \cite[Theorem 2.54]{ceccherini-silberstein-scarabotti-tolli2022:RepresentationFiniteGroupExtensions}.
\begin{lemma}\label{lem:little-groups}
  Let \(G=A\rtimes H\) be a semidirect product of a group \(H\) with an abelian group \(A\) and
  suppose \((V, \rho)\) is an irreducible representation of \(A\).
  \begin{thmlist}
    \item The inertia subgroup \(I_{G}([V])\) is a semidirect product \(A \rtimes H_{[V]}\) with \(H_{[V]}\subset H\).
    \item If \(W \in \widehat{I_{G}([V])}_{[V]}\), there exists an irreducible representation \((U, \tau)\) of \(H_{[V]}\) such that \(W\cong (V\otimes U, \varphi)\), where
    \begin{equation*}
      \varphi \from A\rtimes H_{[V]} \to \Aut_{ \Bbbk}(V\otimes U), \qquad \varphi((a,h))(v\otimes u) = \rho(a)v\otimes \tau(h)u.
    \end{equation*}

  \end{thmlist}
\end{lemma}

Moreover, we will need the following observation about irreducible representations of products of groups.

\begin{lemma}\label{lem:prod-group-trep-reps}
  Let \(G= N \times H\).
  Every irreducible representation \(U\) of \(G\) is isomorphic to a tensor product \(V\otimes W\) of irreducible representations \(V \in \Rep{N}\) and \(W\in \Rep{H}\).
\end{lemma}
\begin{proof}
  As \(nh =hn\) for all \(n \in N\) and \(h\in H\) the inertia subgroup of any  \([V]\in \widehat{N}\) is \(I_{G}([V]) = G\) and the claim follows from \cite[Theorem 2.26]{ceccherini-silberstein-scarabotti-tolli2022:RepresentationFiniteGroupExtensions}. Alternatively, one can argue that since \( \Bbbk\) is algebraically closed and $\Bbbk[G] \simeq \Bbbk[N] \otimes \Bbbk[H]$ is semisimple, the Artin-Wedderburn decomposition of $\Bbbk[G]$ is obtain by the direct sum of the tensor products of the Artin-Wedderburn blocks of $\Bbbk[N]$ with the blocks of $\Bbbk[H]$.
\end{proof}

\subsection{Irreducible representations of symmetric groups}\label{sec:simple-representations-of-symmetric-groups}
We will show in Section~\ref{sec:proof-of-the-main-theorem} that the inertia subgroups arising in the study of the generalised Kac--Paljutkin Hopf algebras are closely related to symmetric groups \(S_{k}\).
In order to determine \(\widehat{S_{k}}\),  we provide a short summary of the theory of Young symmetrizers based on \cite{howe2022:InvitationToRepresentationTheory}.

\begin{definition}\label{def:partition}
  Consider a sequence \(\mu = (\mu_{i})_{i \in \mathbb{N}}\) of non-increasing non-negative integers
  In case \(|\mu| \eqdef \Sigma_{i\in \mathbb{N}} \mu_{i} = k\), we call \(\mu\) a \emph{partition of \(k\)} and we write \(\mu \vdash k\).

  We refer to \(\ell(\mu)= \max\{i \in \mathbb{N} \mid \mu_{i} \neq 0\}\) as the \emph{length} of \(\mu\).
\end{definition}
We can visualise each partition \(\mu \vdash k\) via its \emph{Young diagram}. That is, a left-justified array of unit-sized boxes with \(\mu_{i}\) boxes in the \(i\)-th row.
For example, the diagram corresponding to  \((3,2,2)\vdash 7\) is:
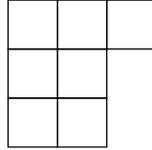
\begin{equation*}
  \ydiagram{3, 2, 2}
\end{equation*}

\begin{definition}\label{def:Young-tableau}
  A \emph{Young tableau} \(T\) of shape \( \mu \vdash k\) is a bijection
  \begin{equation*}
    T\from \{(i,j) \in \mathbb{N}^{2}\mid j\leq \mu_{j}\} \to [k], \qquad (i,j) \mapsto T_{i,j}.
  \end{equation*}
  It is called \emph{standard} if for all \((i,j)\in \mathbb{N}^{2}\) we have \(T_{i,j}<T_{i,j+1}\) and \(T_{i,j}< T_{i+1,j}\).
\end{definition}

A Young tableau of shape \(\mu \vdash k\) can be graphically interpreted as a Young diagram whose boxes are filled with the numbers \(1, \dots, k\).
The leftmost of the three diagrams shown below is in standard form, while in the middle one the rows, and in the rightmost one the columns are not sorted in ascending order.
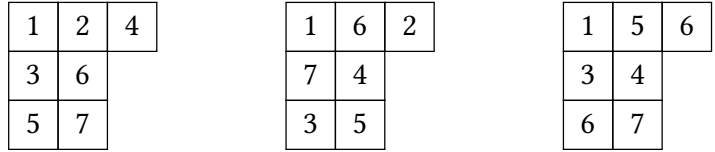
\begin{equation*}
  \ytableaushort{ 1 2 4, 3 6, 5 7} \qquad\qquad
  \ytableaushort{ 1 6 2, 7 4, 3 5} \qquad\qquad
  \ytableaushort{ 1 5 6, 3 4, 6 7}
\end{equation*}

\begin{definition}\label{def:hor-vert-subgroups}
  Consider a Young tableau \(T\) of shape \(\mu \vdash k\).
  Its \emph{horizontal group} and \emph{vertical group}  are respectively given by
  \begin{equation*}
    H_{T} \eqdef \{\sigma \in S_{k} \mid \sigma(T_{i,j}) = T_{i, j'} \} \quad \text{ and }\quad
    V_{T} \eqdef \{\sigma \in S_{k} \mid \sigma(T_{i,j}) = T_{i',j} \}.
  \end{equation*}
\end{definition}

By definition, the horizontal group of a tableau \(T\) preserves the rows of its Young diagram and the vertical group its columns.

\begin{definition}\label{def:Young-symmetriser}
  Let \(T\) be a Young tabelau of shape \(\mu \vdash k\) and \(f_{\mu}\in \mathbb{N}\) the number of standard Young tableaux of shape \(\mu\).
  The \emph{(normalised) Young symmetrizer} of \(T\) is
  \begin{equation*}
    e_{T}\eqdef \tfrac{f_{\mu}}{k!} h_{T} v_{T}\in \Bbbk [S_k], \qquad \text{where }
    h_{T} = \sum_{g\in H_{T}} g \quad \text{ and }\quad
    v_{T} =  \sum_{g\in V_{T}} \mathrm{sgn}(g)g.
  \end{equation*}
\end{definition}

The following result is proven in Proposition~8.30, Theorem~8.22, and Corollary~11.21 of \cite{howe2022:InvitationToRepresentationTheory}.
\begin{lemma}\label{lem:specht-modules}
  Let \(T\) be a Young tableau of shape \(\mu \vdash k\). Then
  \begin{thmlist}
    \item \(e_{T}\) is an idempotent,
    \item \(V^{\mu}\eqdef \mathbb{k}[S_{k}]e_{T}\) is an irreducible  representation of \(S_{k}\),
    \item  we have \(\dim V^{\mu} =|\{\text{standard tableaux of shape } \mu\}|=f_\mu\), and
    \item if \(R\) is another tableau of shape \(\mu'\), we have \(\mathbb{C}[S_{k}]e_{T} \cong \mathbb{C}[S_{k}]e_{R}\) if and only if \(\mu = \mu'\).
  \end{thmlist}
  Moreover, any irreducible representation \(W\) of \(S_{k}\) is isomorphic to some  \(V^{\mu'}\) with \(\mu' \vdash k\).
\end{lemma}

\section{Proof of the main theorem}\label{sec:proof-of-the-main-theorem}
Let \(\mathcal{Y}_{k}\) be the set of all partitions of \(k\), define \(\mathcal{Y}_{0}=\{\emptyset\}\), and write \(\mathcal{Y}\eqdef \cup_{k\in \mathbb{N}_{0}} \mathcal{Y}_{k}\).
\begin{definition}\label{def:multipartition}
  We call a map \(\beta\from [n] \to \mathcal{Y} \), \(i\mapsto \mu_{i} \) such that \(\sum_{i\in [n]} |\mu_{i}|=m\) holds an \emph{\( [n]\)-labelled partition of \(m\)}.
\end{definition}

  Consider the  \([3]\)-labelled partition \(\beta\from [3] \to \mathcal{Y} \) of \(10\) given by
  \begin{equation*}
    1\mapsto (3,2,2) \vdash 7, \qquad
    2 \mapsto \emptyset, \qquad
    3 \mapsto (1,1,1)\vdash 3.
  \end{equation*}
  In order to construct an irreducible \(H_{3,10}\)-module  \(V^{\beta}\) of out of \(\beta\), we will now fill the boxes of its diagrams as shown below and use that to construct an idempotent \(e_{\beta}\in H_{3,10}\):
  \begin{equation*}
    \begin{aligned}
      \ytableaushort{ \none {\none[0]} \none} & \quad\ytableaushort{ \none {\none[1]} \none}\quad && \ytableaushort{ \none {\none[2]} \none}\\
      \ytableaushort{1 2 3, 4 5, 6 7}  & \quad \ytableaushort{ \none {\none[*]} \none} && \ytableaushort{\none 1, \none {2}, \none {3}}
    \end{aligned}
  \end{equation*}

\begin{construction}\label{const:idempotent-from-beta}
  Let \(\beta\) be an \([n]\)-labelled partition of \(m\) and write \(\mu_{i}= \beta(i)\) and \(l_i = |\mu_i|\).
  \begin{thmlist}
    \item We associate to \(\beta\) the element \(\mathbf{\lambda}_{\beta} =\mathbf{\lambda} \in \mathbb{Z}_{n}^{m}\) defined by
    \[ \mathbf{\lambda}_j \eqdef \max\left\{ i : \sum_{k=1}^i l_k < j\right\}, \qquad \mbox{ for  } j \in [n], \mbox{i.e. }\]
  \begin{equation*}
   \mathbf{\lambda}_{1} = \dots = \mathbf{\lambda}_{l_1} = 0, \quad
    \mathbf{\lambda}_{l_1+1} = \dots  = \mathbf{\lambda}_{l_1+l_2}= 1, \quad\dots,\quad
    \mathbf{\lambda}_{m-l_{n}+1} = \dots = \mathbf{\lambda}_{m} = n - 1.
  \end{equation*}
  \item We denote by \(S_{\beta}\subset H_{n,m}\) the subalgebra generated by the elements
  \begin{equation*}
    s_{1}, \dots , s_{l_1-1},\quad
    s_{l_1+1}, \dots , s_{l_1+l_2-1},\quad
    \dots ,\quad
    s_{m-l_n+1}, \dots , s_{m-1}.
  \end{equation*}
  It is isomorphic to \(\Bbbk [S_{l_1}\times \dots \times S_{l_n}]\) and we write \(\iota_{i}\from  \Bbbk[S_{l_{i}}] \to S_{\beta}\) for the canonical inclusion of \( \Bbbk[S_{l_i}]\) into  \(S_{\beta}\).
  \item We choose for each \(i \in [n]\) the unique standard Young tableau \(T^{(i)}\) of  shape \(\mu_{i}\) satisfying \(T^{(i)}_{(j,k)}+1=T^{(i)}_{(j,k+1)}\).
  \item We set
  \begin{equation}\label{eq:the-idempotent}
    e_{\beta}\eqdef \Lambda_{\mathbf{\lambda}_{\beta}} \iota_{1}(e_{T^{(1)}}) \cdot \dots \cdot \iota_{n}(e_{T^{(n)}}).
  \end{equation}
  \end{thmlist}
\end{construction}

\begin{theorem}\label{thm:the-one-we-want}
  The map
  \begin{equation}
    \begin{aligned}
      \{[n]\text{-labelled partitions of }m\}& \to \widehat{H_{n,m}}\\
      \beta &\mapsto \big[H_{n,m} e_{\beta}\big]
    \end{aligned}
  \end{equation}
  is a bijection.
\end{theorem}
\begin{proof}
  To keep our arguments concise, we write \(N\) and \(S\) for the subgroups of \(H_{n,m}\) generated by \(x_{1}, \dots, x_{m}\) and \(s_{1}, \dots , s_{m-1}\), respectively.
  Furthermore, we define \(G= N\wr S\).

  Using Proposition~\ref{prop:KP-is-hyperoct}, we obtain
  \begin{equation*}
    \widehat{H_{n,m}} = \widehat{G} = \bigcup_{[V]\in \widehat{N}} \widehat{G}_{[V]}.
  \end{equation*}
  Suppose \(V, V'\) are two irreducible representations of \(N\).
  Lemma~\ref{lem:Cliffords-thm}  implies that \(\widehat{G}_{[V]}\) and \(\widehat{G}_{[V']}\) coincide if there exists a \(g\in G\) such that \(V' \cong {}^{g}V\).
  Otherwise, \(\widehat{G}_{[V]}\) and \(\widehat{G}_{[V']}\) are disjoint.

  We identify \(\widehat{N}\) with \( \mathbb{Z}_{n}^{m}\) using the bijection
  \begin{equation*}
    \mathbb{Z}_{n}^{m} \to \widehat{N},\qquad \mathbf{\lambda} \mapsto[ \Bbbk[N] \Lambda_{\mathbf{\lambda}}].
  \end{equation*}
  Equation~\eqref{eq:z-action-on-simples} shows that \(G\) acts transitively on \(\widehat{N}\).
  Therefore, each \(G\)-orbit of \(\widehat{N}\) contains a unique element \(\mathbf{\lambda} = (\lambda_{1}, \dots , \lambda_{m})\), where \(\lambda_{1} \leq \dots \leq \lambda_{m}\) is a non-decreasing tuple of elements \(\lambda_{i}\in [n]\).
  For each \(i \in [n]\), we write \(k_{i} = |\{\lambda_j \in \mathbf{\lambda} \mid \lambda_{j} = i\}|\).
  A direct computation shows that the inertia subgroup \(I_{G}(\mathbf{\lambda})\) is generated by
  \begin{gather*}
    x_{1}, \dots , x_{m}, \\
    s_{1}, \dots , s_{k_{1}-1}, \quad  s_{k_{1}+1}, \dots , s_{k_{1}+k_{2}-1},\quad  \dots ,\quad  s_{m- k_{n}+1}, \dots , s_{m-1}
  \end{gather*}
  and we can canonically identify it with the semidirect product
  \(\mathbb{Z}_{n}^{m}\rtimes (S_{k_{1}} \times \dots \times S_{k_{n}})\).
  By combining~Lemmas~\ref{lem:little-groups} and~\ref{lem:prod-group-trep-reps}, we observe that any \(U\in\Rep{I_{G}(\mathbf{\lambda})}\) with  \([U]\in (\widehat{I_{G}(\mathbf{\lambda})})_{\mathbf{\lambda}}\) is isomorphic to a tensor product \(  \Bbbk[N] \Lambda _{\mathbf{\lambda}}\otimes V_{1} \otimes \dots \otimes V_{n}\) with \(V_{i}\) an irreducible representation of \(S_{k_{i}}\).
  Lemma~\ref{lem:specht-modules} shows that there exists a unique Young diagram \(\mu_{i} \vdash k_{i}\) such that \(V_{i} =  \Bbbk[S_{k_{i}}] e_{\mu_{i}}\).
  We write \(\beta\from [n] \to \mathcal{Y}\) for the \([n]\)-labelled partition defined by \(\beta(i) = \mu_{i}\).
  By Equation~\eqref{eq:z-action-on-simples}, we have for each \(s \in  \Bbbk[S_{k_{i}}]\), \(t\in  \Bbbk[S_{k_{j}}]\) that \(s \Lambda_{\mathbf{\lambda}} = \Lambda_{\mathbf{\lambda}} s\) and \(st = ts\) if \(i \neq j\).
  Fixing for each \(1\leq i \leq n\) the standard Young tableau \(T^{(i)}\) of shape \(\mu_{i}\) such that
  \(T^{(i)}_{a,b}+1=T^{(i)}_{a, b+1}\), we get
  \begin{align*}
    U & \cong \Bbbk[ N]\Lambda_{\mathbf{\lambda}} \otimes  \Bbbk[S_{k_{1}}]e_{T^{(1)}}\otimes \dots \otimes  \Bbbk[S_{k_{n}}]e_{T^{(n)}}\\
        &\cong \Bbbk[N\rtimes (S_{k_{1}} \times \dots \times S_{k_{n}})]  \Lambda_{\mathbf{\lambda}} \iota_{1}(e_{T^{(1)}}) \dots \iota_{n}(e_{T^{(n)}}) \\
      & = \Bbbk[I_{G}(\mathbf{\lambda})]  \Lambda_{\mathbf{\lambda}} \iota_{1}(e_{T^{(1)}}) \dots \iota_{n}(e_{T^{(n)}})\\
       & = \Bbbk[I_{G}(\mathbf{\lambda})] e_{\beta}
  \end{align*}

  Finally, since \(e_{\beta}\) is an idempotent, we get
  \begin{equation*}
    \Ind_{I_{G}(\lambda)}^{H_{n,m}} U \cong H_{n,m} \otimes_{\Bbbk[I_{G}(\mathbf{\lambda})]} \Bbbk[I_{G}(\mathbf{\lambda})] e_{\beta} \cong H_{n,m} e_{\beta}
  \end{equation*}
  and the claim follows by Lemma~\ref{lem:clifford-correspondence}.
\end{proof}
As applications of the previous result, we determine the number of isomorphism classes of irreducible representations and state the dimensions of these representations.

\begin{corollary}\label{cor:num-of-irreds}
  Let $m, n \geq 2$.
  The number of irreducible representations of \(H_{n,m}\) is equal to the number of conjugacy classes of the generalised symmetric group \(\mathbb{Z}_n \wr S_m\) and equal to
  \begin{equation}\label{eq:number-of-simples}
    \sum_{ (l_1,\dots, l_n) } p(l_1) \cdots p(l_n)
  \end{equation}
  where the sum runs over all $n$-tuples $(l_1,\dots, l_n)$ with $l_1+\dots+l_n=m$ and where  $p(-)$ denotes the partition function.
\end{corollary}
\begin{proof}
The number of $[n]$-labelled partitions of $m$, i.e. functions $\beta:[n]\to \mathcal{Y}$, such that
$\beta(i) \vdash l_i \eqdef |\beta(i)|$ is a partition of $l_i$ and $l_1+\cdots +l_n = m$, is given by
\(  \sum_{ (l_1, \dots, l_n) } p(l_1) \cdots p(l_n) \), where the sum runs over all $n$-tuples $(l_1,\dots, l_n)$ with $l_1+\cdots+l_n=m$.
Theorem \ref{thm:the-one-we-want} shows that there is a bijection between $[n]$-labelled partitions of $m$ and irreducible representations of $H_{n,m}$. Hence, this expression calculates the number of non-isomorphic irreducible representations, which is also equal to  the number of conjugacy classes of  \(\mathbb{Z}_n \wr S_m\), see \eg \cite[Theorem~3.8]{Kerber}.
\end{proof}

\begin{corollary}\label{cor:dimension-of-simple-modules}
  Let \(\beta\from[n] \to \mathcal{Y}\) be an \([n]\)-labelled partition. We  write \(l_{i} \eqdef |\beta(i)|\) as well as \(f_{i}\eqdef |\{ \text{standard Young tabeaux of shape } \beta(i)\}|\) for all \(i \in [n]\).
  The irreducible representation \(H_{n,m} e_{\beta}\) has the dimension
  \begin{equation}
    \dim H_{n,m} e_{\beta} = \frac{m!}{l_{1}!\cdot\dots \cdot l_{n}!}f_{1} \cdot \dots \cdot f_{n}.
  \end{equation}
\end{corollary}
\begin{proof}
  In line with Construction~\ref{const:idempotent-from-beta}, we associate to \(\beta\) the element \(\mathbf{\lambda}_{\beta}\in \mathbb{Z}_{n}^{m}\) defined by
  \begin{equation*}
    \lambda_{1} = \dots = \lambda_{l_{1}} = 0, \quad
    \lambda_{l_{1}+1} = \dots =\lambda_{l_{1}+l_{2}}= 1, \quad\dots,\quad
    \lambda_{m-l_{n}+1} = \dots = \lambda_{m} = n - 1.
  \end{equation*}
  Its inertia subgroup \(I_{G}(\mathbf{\lambda})\) is \( \mathbb{Z}_{n}^{m}\rtimes S_{l_{1}}\times \dots \times S_{l_{n}}\).
  We fix for every \(i\in [n]\) a standard Young tableau \(T_{i}\) of shape \(\beta(i) \vdash l_{i}\) and write \(V_{i}\eqdef  \Bbbk[S_{i}] e_{T_{i}}\).
  Note that \(\Bbbk[I_{G}(\mathbf{\lambda})] e_{\beta} \cong V_{1}\otimes \dots \otimes V_{n}\) and Lemma~\ref{lem:specht-modules} implies that
  \begin{equation*}
    \dim \Bbbk[I_{G}(\mathbf{\lambda})] e_{\beta}
    = \dim V_{1}\otimes \dots \otimes V_{n}
    = f_{1} \cdot \dots \cdot f_{n}.
  \end{equation*}
  Since \(H_{n,m} e_{\beta} \cong H_{n,m} \otimes_{\Bbbk[I_{G}(\mathbf{\lambda})]} \Bbbk[I_{G}(\mathbf{\lambda})] e_{\beta}\) and \(H_{n,m}\cong  \Bbbk[ \mathbb{Z}_{n} \wr S_{m}]\) is a free \( \Bbbk[I_{G}(\mathbf{\lambda})]\)-module of rank
  \begin{equation*}
    \frac{\dim H_{n,m}}{\dim  \Bbbk[I_{G}(\mathbf{\lambda})]}
    = \frac{n^{m}m!}{n^{m}(l_{1}! \cdot \dots \cdot l_{n}!)}
    = \frac{m!}{l_{1}! \cdot \dots \cdot l_{n}!},
  \end{equation*}
  we get
  \begin{equation*}
    \dim H_{n,m} e_{\beta} = \frac{m!}{l_{1}! \cdot \dots \cdot l_{n}!} f_{1} \cdot \dots \cdot f_{n}. \qedhere
  \end{equation*}
\end{proof}

\begin{remark}
Let $\lambda \vdash l$ be a partition.
Given a box $(a,b)$ in a Young diagram of $\lambda$, a hook $H_{a,b}$ is the box $(a,b)$, along with the boxes in its row to the right and the boxes in its column below. The number of boxes in a hook, $h(a,b) \eqdef | H_{a,b} |$, is called the hook length. The Hook Length formula (\cite[Proposition 13.4]{howe2022:InvitationToRepresentationTheory}) says that the number of standard Young tableaux with shape $\lambda$ is equal to $l!$ divided by the product of the hook lengths $h(a,b)$ of each box $(a,b)$ in the Young diagram of $\lambda$. Hence, Corollary \ref{cor:dimension-of-simple-modules} says that for
 \(\beta\from[n] \to \mathcal{Y}\), an $[n]$-labelled partition of $m$ one has
  \begin{equation*}
    \dim H_{n,m} e_{\beta} = \frac{m!}{ \prod_{i\in[n]}\prod_{(a,b) \in \mathbb{N}^{2}: b\leq \beta(i)_{a}} h(a, b)}.
  \end{equation*}

\end{remark}

To illustrate our construction, we consider the \(48\)-dimensional algebra \(H_{2,3}\) and explicitly determine all of its irreducible representations.

\begin{example}\label{ex:n=2,m=3}
  By Corollary~\ref{cor:num-of-irreds}, there are \(10\) isomorphism classes of irreducible representations of \(H_{2,3}\).
  Explicitly, they correspond to the \([2]\)-partitions of \(3\) shown below.
  
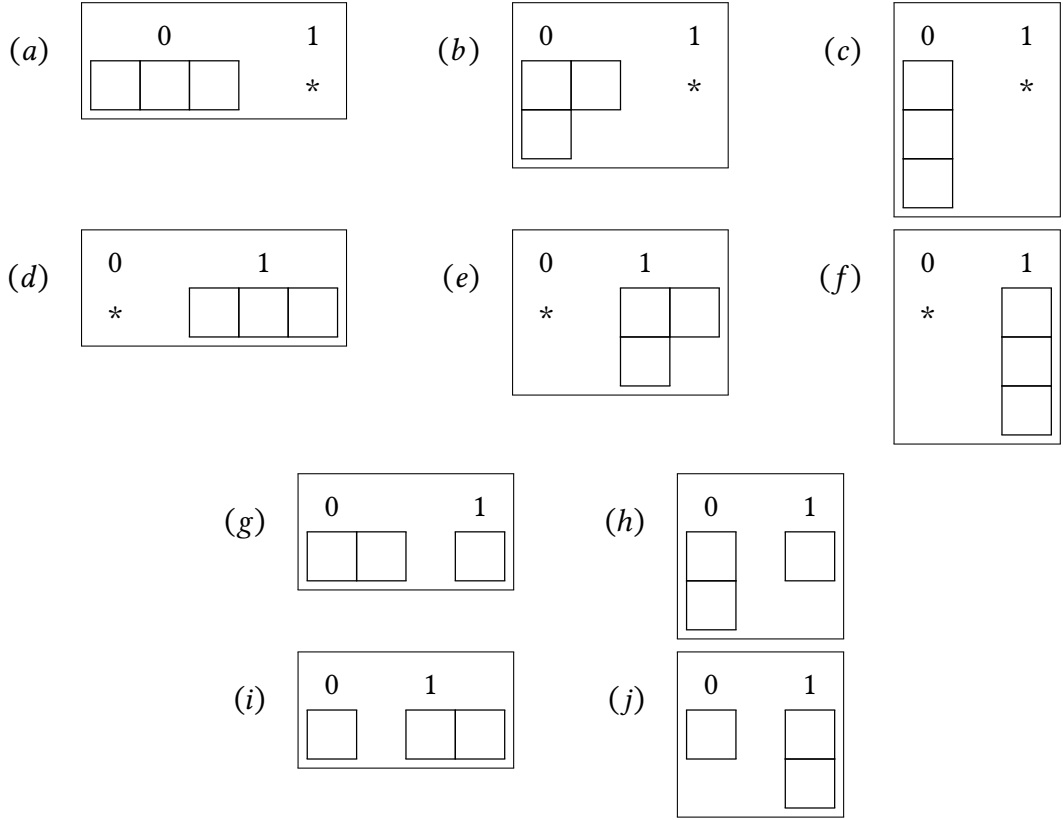
\begin{gather*}
    \begin{aligned}
      &&(a) \quad &\boxed{\ytableaushort{{\none} {\none[0]} {\none} {\none} {\none[1]}, { } { } {} {\none} {\none[\ast]}}}\qquad
      &(b) \quad &\boxed{\ytableaushort{{\none[0]} {\none} {\none} {\none[1]}, { } { } {\none} {\none[\ast]} , { }}}\qquad
      &(c) \quad &\boxed{\ytableaushort{{\none[0]} {\none} {\none[1]}, { } {\none} {\none[\ast]} , { }  , { }}}\\
      &&(d) \quad &\boxed{\ytableaushort{{\none[0]} {\none} {\none} {\none[1]} {\none},{\none[*]} {\none} { } { }{}}}
      &(e) \quad &\boxed{\ytableaushort{{\none[0]} {\none} {\none[1]}{\none},{\none[\ast]} {\none} { } { }, {\none} {\none} { }}}
      &(f) \quad &\boxed{\ytableaushort{{\none[0]} {\none} {\none[1]}, {\none[\ast]} {\none} { },{\none}{\none} { }  , {\none} {\none} { }}}
    \end{aligned}
  \end{gather*}
  \begin{gather*}
    \begin{aligned}
      &&(g) \quad &\boxed{\ytableaushort{{\none[0]} {\none} {\none} {\none[1]}, { } {} {\none} {}}}\qquad
      &(h) \quad & \boxed{\ytableaushort{{\none[0]} {\none} {\none[1]}, { } {\none} {} , { }}}\\
      &&(i) \quad &\boxed{\ytableaushort{{\none[0]} {\none} {\none[1]} {\none}, { } {\none} {} {}}}\qquad
      &(j) \quad &\boxed{\ytableaushort{{\none[0]} {\none} {\none[1]}, { } {\none} {} , {\none} {\none} { }}}
    \end{aligned}
  \end{gather*}
  Construction~\ref{const:idempotent-from-beta} yields the corresponding idempotents and Corollary~\ref{cor:dimension-of-simple-modules} allows us to compute the dimensions of the associated irreducible representations:
  \[\begin{array}{lllc}
    \text{number} & \beta            &   e_\beta                                             & \text{dimension} \\\hline
    (a)           & ((3),*)      & \frac{1}{6}\Lambda_{(0,0,0)}(1+s_1+s_2+s_1s_2+s_2s_1+s_1s_2s_1) & 1 \\[1mm]
    (b)           & ((2,1),*)    & \frac{1}{3}\Lambda_{(0,0,0)}(1+s_1)(1-s_1s_2s_1)             & 2 \\[1mm]
    (c)           & ((1,1,1),*)  & \frac{1}{6}\Lambda_{(0,0,0)}(1-s_1-s_2+s_1s_2+s_2s_1-s_1s_2s_1) & 1 \\[1mm]
    (g)           & ((2),(1))    & \frac{1}{2}\Lambda_{(0,0,1)}(1+s_1)                        & 3\\[1mm]
    (h)           & ((1,1),(1))  & \frac{1}{2}\Lambda_{(0,0,1)}(1-s_1)                        & 3\\[1mm]\hline
    (d)           & (*, (3))     & \frac{1}{6}\Lambda_{(1,1,1)}(1+s_1+s_2+s_1s_2+s_2s_1+s_1s_2s_1) & 1\\[1mm]
    (e)           & (*, (2,1))   & \frac{1}{3}\Lambda_{(1,1,1)}(1+s_1)(1-s_1s_2s_1)             & 2\\[1mm]
    (f)           & (*, (1,1,1)) & \frac{1}{6}\Lambda_{(1,1,1)}(1-s_1-s_2+s_1s_2+s_2s_1-s_1s_2s_1) & 1\\[1mm]
    (i)           & ((1),(2))    & \frac{1}{2}\Lambda_{(1,1,0)}(1+s_2)                        & 3\\[1mm]
    (j)           & ((1),(1,1))  & \frac{1}{2}\Lambda_{(1,1,0)}(1-s_2)                        & 3
  \end{array}\]
Note that the horizontal line emphasises the symmetry of this table obtained by tensoring with the one-dimensional simple representation corresponding to the \([2]\)-labelled partition of \(3\) shown in diagram \((d)\).
\end{example}

\printbibliography

\end{document}